\documentclass[12pt]{article}                                               

\input{amssym.def}                                                           
\input{amssym.tex}

\begin{document}                                                             
\title{Real multiplication revisited}

\author{Igor  ~V. ~Nikolaev
\footnote{Partially supported 
by NSERC.}
}


\date{}
 \maketitle


\newtheorem{thm}{Theorem}
\newtheorem{lem}{Lemma}
\newtheorem{dfn}{Definition}
\newtheorem{rmk}{Remark}
\newtheorem{cor}{Corollary}
\newtheorem{cnj}{Conjecture}
\newtheorem{exm}{Example}

\newcommand{\Qcoh}{\hbox{\bf Qcoh}}

\begin{abstract}
It is proved that  the Hilbert class field of a 
real quadratic field ${\Bbb Q}(\sqrt{D})$ modulo a power $m$ of the conductor $f$ 
is generated by the Fourier coefficients of the  Hecke eigenform  for a congruence subgroup of  level $fD$.

\vspace{7mm}


\vspace{5mm}

\end{abstract}

\section{Introduction}
The  Kronecker's {\it Jugendtraum} is a conjecture that the maximal unramified abelian extension
(the Hilbert class field)  of any algebraic number
field is generated by the special values of modular functions attached  to an 
abelian variety.   The conjecture is true for the rational   field  and  
imaginary quadratic fields   with  the modular functions 
being an exponent    and the $j$-invariant,  respectively.  
 In the case of an arbitrary number field,  a description of the abelian extensions is given 
 by class field theory, but an explicit formula for  the generators of these abelian extensions, 
 in the sense sought by Kronecker,   is unknown even for the real quadratic fields.

The problem was first studied by   [Hecke 1910]   \cite{Hec1}.  A description of abelian extensions of real 
quadratic number fields in terms of coordinates of points of finite order on abelian varieties associated 
with certain modular curves was obtained in [Shimura  1972]   \cite{Shi1}.  
Stark  formulated a number of conjectures on abelian extension of arbitrary number fields,  
 which in the real quadratic case amount to specifying generators of these extensions using special values of 
 Artin $L$-functions,   see  [Stark 1976]   \cite{Sta1}.  
Based on an analogy with complex multiplication,  Manin  suggested to use  the so-called  
``pseudo-lattices''  ${\Bbb Z}+{\Bbb Z}\theta$ in  ${\Bbb R}$   having  non-trivial real multiplications
to produce abelian extensions  of real quadratic fields,    see [Manin 2004]   \cite{Man1}.
 Similar to the case of complex multiplication,  the endomorphism ring ${\goth R}_{\goth f}={\Bbb Z}+{\goth f} O_{\goth k}$  
 of pseudo-lattice  ${\Bbb Z}+{\Bbb Z}\theta$  is an order  in the real quadratic field ${\goth k}={\Bbb Q}(\theta)$,  
 where  $O_{\goth k}$  is the ring of integers of ${\goth k}$  and ${\goth f}$  is the conductor of ${\goth R}_{\goth f}$;
Manin calls these pseudo-lattices with {\it real multiplication}.

The aim of our note is a  formula  for  generators of  the Hilbert class field   of 
 real quadratic fields  based on a modularity and a  symmetry 
 of complex and real multiplication.     To give an idea,  let 
\begin{equation}
\Gamma_1(N)=\left\{\left(\matrix{a & b\cr c & d}\right)\in SL_2({\Bbb Z}) ~|~
a\equiv d\equiv 1 ~\hbox{{\bf mod}} ~N,  ~c\equiv 0 ~\hbox{{\bf mod}}  ~N\right\}
\end{equation}
be a congruence subgroup of level $N\ge 1$ and  ${\Bbb H}$ 
be  the Lobachevsky half-plane;  let  $X_1(N):={\Bbb H}/\Gamma_1(N)$ be the corresponding 
modular curve and $S_2(\Gamma_1(N))$ the space of all cusp forms on $\Gamma_1(N)$ of weight 2.
 Let  ${\cal E}_{CM}^{(-D,f)}$  be elliptic curve with complex multiplication
by an order $R_f={\Bbb Z}+fO_k$ in  the  field $k={\Bbb Q}(\sqrt{-D})$,  see 
[Silverman 1994] \cite{S},  Chapter II.   Denote by ${\cal K}^{ab}(k):=k(j({\cal E}_{CM}^{(-D,f)}))$ 
the Hilbert class field  of  $k$ modulo conductor   $f\ge 1$ and let  $N=fD$; 
let   $Jac~(X_1(fD))$ be  the Jacobian of modular curve $X_1(fD)$. 
There exists an abelian sub-variety $A_{\phi}\subset  Jac~(X_1(fD))$,
such that its points of finite order generate ${\cal K}^{ab}(k)$, 
see  [Hecke 1928]  \cite{Hec2},  [Shimura 1971] \cite{Shi2},  Theorem 1 and [Shimura 1972]  \cite{Shi1},   
Section 8.   The  ${\cal K}^{ab}(k)$ is a {\it CM-field},  i.e. a totally imaginary quadratic extension of the totally real field 
${\cal K}_{\phi}$   generated by the Fourier coefficients of the Hecke eigenform 
$\phi(z)\in S_2(\Gamma_1(fD))$,   see [Shimura 1972]  \cite{Shi1}, p. 137. 
In particular,    there exists a holomorphic map $X_1^0(fD)\to {\cal E}_{CM}^{(-D,f)}$,
where $X_1^0(fD)$ is a Riemann surface such that $Jac~(X_1^0(fD))\cong A_{\phi}$;
we refer to the above as a {\it modularity} of complex multiplication.

Recall that  (twisted homogeneous)  coordinate ring of an elliptic curve ${\cal E}({\Bbb C})$ is isomorphic to a 
{\it Sklyanin algebra},   see  e.g.  [Stafford \& van ~den ~Bergh  2001]  \cite{StaVdb1},
Example 8.5;  the norm-closure of a self-adjoint representation of the Sklyanin algebra
 by the linear operators on a Hilbert space ${\cal H}$  is isomorphic to a 
noncommutative torus ${\cal A}_{\theta}$,  see [Rieffel  1990]  \cite{Rie1} for the definition.    
Whenever elliptic curve ${\cal E}({\Bbb C})\cong  {\cal E}_{CM}^{(-D,f)}$ has complex multiplication,  
the noncommutative torus ${\cal A}_{\theta}$  has real multiplication by an order 
${\goth R}_{\goth f}={\Bbb Z}+{\goth f}O_{\goth k}$   in the field ${\goth k}={\Bbb Q}(\sqrt{D})$;
moreover, it is known that ${\goth f}=f^m$ for the minimal power $m$ satisfying   an isomorphism:
\begin{equation}\label{eq3}
Cl~({\goth R}_{f^m})\cong Cl~(R_f),
\end{equation}
where  $Cl~(R_f)$ and  $Cl~({\goth R}_{\goth f})$ are the ideal  class groups of orders $R_f$ and 
${\goth R}_{\goth f}$, respectively.    We shall refer to (\ref{eq3})  as a {\it symmetry} of complex and real multiplication.
 The noncommutative torus with real multiplication
by  ${\goth R}_{\goth f}$ will be denoted by ${\cal A}_{RM}^{(D, {\goth f})}$. 
\begin{rmk}
\textnormal{
The isomorphism (\ref{eq3}) can be calculated  using  
the well-known  formula for the class number  of a non-maximal  order ${\Bbb Z}+fO_K$ of a quadratic
field $K={\Bbb Q}(\sqrt{D})$:
\begin{equation}
h_{{\Bbb Z}+fO_K}={h_{O_K} f \over e_f}\prod_{p|f}\left(1-\left({D\over p}\right){1\over p}\right),
\end{equation}
where $h_{O_K}$ is the class number of the maximal order $O_K$,  $e_f$ is the index
of the group of units of  ${\Bbb Z}+fO_K$ in the group of units of $O_K$, $p$ is a prime number and 
$\left({D\over p}\right)$ is the Legendre symbol,   see e.g.  [Borevich \& Shafarevich 1966]  
\cite{BS}, p.153  and [Hasse 1950]  \cite{H},  pp. 297 and 351.   
}
\end{rmk}
The (twisted homogeneous) coordinate ring of the Riemann surface $X_1^0(fD)$ is  
an  {\it AF-algebra}   ${\Bbb A}_{\phi^0}$ linked  to a holomorphic differential $\phi^0(z)dz$
on $X_1^0(fD)$,  see Section 2.2,   Definition \ref{dfn1} and Remark \ref{rmk4} for the details;
the  Grothendieck semigroup $K_0^+({\Bbb A}_{\phi^0})$ is a pseudo-lattice  
${\Bbb Z}+{\Bbb Z}\theta_1+\dots+{\Bbb Z}\theta_{n-1}$ in the number field  
 ${\cal K}_{\phi}$,  where    $n$ equals  the  genus of $X_1^0(fD)$. 
Moreover,  a holomorphic  map $X_1^0(fD)\to {\cal E}_{CM}^{(-D,f)}$
 induces the  $C^*$-algebra homomorphism ${\Bbb A}_{\phi^0}\to {\cal A}_{RM}^{(D, {\goth f})}$
 between  the corresponding coordinate rings,  so that  the following diagram commutes:

\begin{picture}(300,110)(-100,-5)
\put(20,70){\vector(0,-1){35}}
\put(130,70){\vector(0,-1){35}}
\put(52,23){\vector(1,0){53}}
\put(52,83){\vector(1,0){53}}
\put(5,20){${\cal E}_{CM}^{(-D,f)}$}
\put(120,20){${\cal A}_{RM}^{(D, {\goth f})}$}
\put(0,80){$X_1^0(fD)$}
\put(127,80){${\Bbb A}_{\phi^0}$}
\put(50, 90){coordinate}
\put(70, 70){map}
\put(50, 30){coordinate}
\put(70, 10){map}
\end{picture}

\noindent
But  $K_0^+({\cal A}_{RM}^{(D, {\goth f})})$ is a pseudo-lattice ${\Bbb Z}+{\Bbb Z}\theta$ in 
the field ${\goth k}$,  such that   $End~({\Bbb Z}+{\Bbb Z}\theta)\cong {\goth R}_{\goth f}$;   
in other words,  
one can   use the above diagram  to control the arithmetic of the field ${\cal K}_{\phi}$ 
by such of  the real quadratic field  $\goth k$. Roughly speaking,
this observation solves the Kronecker's Jugendtraum for the real quadratic fields; namely,
the following  is true. 
\begin{thm}\label{thm1}
The Hilbert class field of a real quadratic field ${\goth k}={\Bbb Q}(\sqrt{D})$ 
modulo conductor $f^m$ is an extension of ${\goth k}$ 
 by the Fourier coefficients  of  the Hecke eigenform   $\phi(z)\in S_2(\Gamma_1(fD))$,
 where $m$ is the smallest  positive integer satisfying  isomorphism (\ref{eq3}). 
\end{thm}
\begin{rmk}
\textnormal{
Theorem \ref{thm1}  can be used to  compute   concrete  extensions.  For  instance,
theorem \ref{thm1}  says  that   for the quadratic field   ${\Bbb Q}(\sqrt{15})$  its  Hilbert class field
is isomorphic to ${\Bbb Q}\left(\sqrt{-1 +\sqrt{15}}\right)$ and 
 for ${\Bbb Q}(\sqrt{14})$   such a  field    modulo conductor ${\goth f}=8$  is isomorphic to  
 ${\Bbb Q}\left(^4\sqrt{-27+8\sqrt{14}}\right)$,   see Section 4 for more examples. 
}
\end{rmk}
The article is organized as follows.   Section 2 covers  basic  facts on real 
multiplication and AF-algebras of the Hecke eigenforms.  Theorem \ref{thm1}
is proved in Section 3.  Section 4 contains  numerical examples illustrating 
theorem \ref{thm1}.

\section{Preliminaries}
The reader can find basics of the $C^*$-algebras in [Murphy 1990]  \cite{M}
and their $K$-theory  in [Blackadar  1986]  \cite{B}.    
The noncommutative tori are covered in [Rieffel  1990]  \cite{Rie1}
and real multiplication in [Manin 2004]  \cite{Man1}.  
For  main ideas  of non-commutative algebraic geometry,  see the survey 
by [Stafford \& van ~den ~Bergh  2001]  \cite{StaVdb1}.  The AF-algebras 
are reviewed in [Effros 1981]  \cite{E}.  For a general theory of modular 
forms we refer the reader to [Diamond \& Shurman 2005]  \cite{DS}.

\subsection{Real multiplication}
The noncommutative torus ${\cal A}_{\theta}$ is a universal {\it $C^*$-algebra}
generated by the unitary operators $u$ and $v$ acting on a Hilbert space ${\cal H}$
and satisfying the commutation relation $vu=e^{2\pi i\theta}uv$,  where $\theta$ is a 
real number.   The $C^*$-algebra ${\cal A}_{\theta}$ is said to be stably isomorphic
(Morita equivalent)  to ${\cal A}_{\theta'}$,  whenever ${\cal A}_{\theta}\otimes {\cal K}\cong 
{\cal A}_{\theta'}\otimes {\cal K}$,  where ${\cal K}$ is the $C^*$-algebra of all compact operators
on ${\cal H}$;  the ${\cal A}_{\theta}$ is stably isomorphic to ${\cal A}_{\theta'}$ if and only if 
\begin{equation}\label{eq6}
\theta'={a\theta +b\over c\theta+d}\quad
\hbox{for some matrix} \quad  \left(\matrix{a & b\cr c & d}\right)\in SL_2({\Bbb Z}). 
\end{equation}
The $K$-theory of ${\cal A}_{\theta}$ is  two-periodic and 
$K_0({\cal A}_{\theta})\cong  K_1({\cal A}_{\theta})\cong {\Bbb Z}^2$ so that  
the Grothendieck  semigroup $K_0^+({\cal A}_{\theta})$ corresponds  to positive reals of 
the pseudo-lattice  ${\Bbb Z}+{\Bbb Z}\theta\subset {\Bbb R}$.
The ${\cal A}_{\theta}$ is said to have {\it real multiplication},    if $\theta$ is a quadratic
irrationality,  i.e.  irrational root of a quadratic polynomial in ${\Bbb Z}[x]$.  
The real multiplication says that the endomorphism ring of pseudo-lattice 
${\Bbb Z}+{\Bbb Z}\theta$ exceeds the ring ${\Bbb Z}$ corresponding to multiplication
by $m$ endomorphisms;  similar to complex multiplication, it means that the
endomorphism ring is isomorphic to an order ${\goth R}_{\goth f}={\Bbb Z}+{\goth f}O_{\goth k}$
of conductor ${\goth f}\ge 1$ in the real quadratic field ${\goth k}={\Bbb Q}(\theta)$,
hence the name.  If $D>0$ is the discriminant of  ${\goth k}$,   then 
by ${\cal A}_{RM}^{(D, {\goth f})}$ we denote torus ${\cal A}_{\theta}$ with real multiplication
by the order ${\goth R}_{\goth f}$.

The Sklyanin algebra $S_{\alpha,\beta,\gamma}({\Bbb C})$  is  a free   ${\Bbb C}$-algebra   on   four  generators  
  and  six   relations: 
\begin{equation}
\left\{
\begin{array}{ccc}
x_1x_2-x_2x_1 &=& \alpha(x_3x_4+x_4x_3),\\
x_1x_2+x_2x_1 &=& x_3x_4-x_4x_3,\\
x_1x_3-x_3x_1 &=& \beta(x_4x_2+x_2x_4),\\
x_1x_3+x_3x_1 &=& x_4x_2-x_2x_4,\\
x_1x_4-x_4x_1 &=& \gamma(x_2x_3+x_3x_2),\\ 
x_1x_4+x_4x_1 &=& x_2x_3-x_3x_2,
\end{array}
\right.
\end{equation}
where $\alpha+\beta+\gamma+\alpha\beta\gamma=0$;  
such an algebra  corresponds to  a  {\it twisted homogeneous  coordinate ring}   of 
an elliptic curve  in the complex projective space ${\Bbb C}P^3$ 
 given by the intersection  of two quadric surfaces of the form
 ${\cal E}_{\alpha,\beta,\gamma}({\Bbb C})=\{(u,v,w,z)\in {\Bbb C}P^3 ~|~u^2+v^2+w^2+z^2={1-\alpha\over 1+\beta}v^2+
 {1+\alpha\over 1-\gamma}w^2+z^2=0\}$.  
 Being such a ring means that   the algebra $S_{\alpha,\beta,\gamma}$ satisfies an isomorphism
\begin{equation}
\hbox{{\bf Mod}}~(S_{\alpha,\beta,\gamma}({\Bbb C}))/
\hbox{{\bf Tors}}\cong \hbox{{\bf Coh}}~({\cal E}_{\alpha,\beta,\gamma}({\Bbb C})),
\end{equation}
 where {\bf Coh} is  the category of quasi-coherent sheaves on ${\cal E}_{\alpha,\beta,\gamma}({\Bbb C})$, 
  {\bf Mod}  the category of graded left modules over the graded ring $S_{\alpha,\beta,\gamma}({\Bbb C})$
 and  {\bf Tors}  the full sub-category of {\bf Mod} consisting of the
torsion modules,  see  [Stafford \& van ~den ~Bergh  2001]  \cite{StaVdb1},  Example 8.5.

If one sets $x_1=u, x_2=u^*, x_3=v, x_4=v^*$,  then there exists a self-adjoint representation
of the Sklyanin $\ast$-algebra $S_{\alpha, 1, -1}({\Bbb C})$ by linear operators on a Hilbert
space ${\cal H}$,  such that its norm-closure is isomorphic to ${\cal A}_{\theta}$; 
 namely,
${\cal A}_{\theta}^0\cong S_{\alpha, 1, -1}({\Bbb C})/I_{\mu}$,   where ${\cal A}_{\theta}^0$ is a dense 
sub-algebra of ${\cal A}_{\theta}$ and $I_{\mu}$ is an ideal generated by the ``scaled  unit''
relations $x_1x_3=x_3x_4={1\over\mu}e$,  where $\mu>0$ is a constant.    Thus the algebra ${\cal A}_{\theta}$ 
is a coordinate ring of elliptic curve ${\cal E}({\Bbb C})$,  such that isomorphic elliptic 
curves correspond to the stably isomorphic (Morita equivalent) noncommutative tori; 
this fact explains  the modular transformation law in (\ref{eq6}).    In particular,  
if ${\cal E}({\Bbb C})$  has complex multiplication  by an  order $R_f={\Bbb Z}+fO_k$ in
a quadratic field $k={\Bbb Q}(\sqrt{-D})$,  then ${\cal A}_{\theta}$ has real multiplication
by an order ${\goth R}_{\goth f}={\Bbb Z}+{\goth f}O_{\goth k}$ in the quadratic field 
${\goth k}={\Bbb Q}(\sqrt{D})$,  where ${\goth f}$ is the smallest  integer satisfying 
an  isomorphism  $Cl~({\goth R}_{\goth f})\cong Cl~(R_f)$, see \cite{Nik2};  
 the  isomorphism   is a necessary and
 sufficient condition for  ${\cal A}_{RM}^{(D, {\goth f})}$ to discern non-isomorphic 
 elliptic curves ${\cal E}_{CM}^{(-D,f)}$ having the same endomorphism ring $R_f$.  
For the constraint ${\goth f}=f^m$, see remark \ref{rmk5}.

\subsection{AF-algebra of the Hecke  eigenform}
An {\it AF-algebra}  (Approximately Finite $C^*$-algebra) is defined to
be the  norm closure of an ascending sequence of   finite dimensional
$C^*$-algebras $M_n$,  where  $M_n$ is the $C^*$-algebra of the $n\times n$ matrices
with entries in ${\Bbb C}$. Here the index $n=(n_1,\dots,n_k)$ represents
the  semi-simple matrix algebra $M_n=M_{n_1}\oplus\dots\oplus M_{n_k}$.
The ascending sequence mentioned above  can be written as 
$M_1\buildrel\rm\varphi_1\over\longrightarrow M_2
   \buildrel\rm\varphi_2\over\longrightarrow\dots,
$
where $M_i$ are the finite dimensional $C^*$-algebras and
$\varphi_i$ the homomorphisms between such algebras.  
The homomorphisms $\varphi_i$ can be arranged into  a graph as follows. 
Let  $M_i=M_{i_1}\oplus\dots\oplus M_{i_k}$ and 
$M_{i'}=M_{i_1'}\oplus\dots\oplus M_{i_k'}$ be 
the semi-simple $C^*$-algebras and $\varphi_i: M_i\to M_{i'}$ the  homomorphism. 
One has  two sets of vertices $V_{i_1},\dots, V_{i_k}$ and $V_{i_1'},\dots, V_{i_k'}$
joined by  $b_{rs}$ edges  whenever the summand $M_{i_r}$ contains $b_{rs}$
copies of the summand $M_{i_s'}$ under the embedding $\varphi_i$. 
As $i$ varies, one obtains an infinite graph called the  {\it Bratteli diagram} of the
AF-algebra.  The matrix $B=(b_{rs})$ is known as  a {\it partial multiplicity matrix};
an infinite sequence of $B_i$ defines a unique AF-algebra.   
An  AF-algebra is called  {\it stationary}  if $B_i=Const=B$,  see [Effros 1981]   \cite{E},  Chapter 6;
when two non-similar matrices $B$ and $B'$ have the same characteristic polynomial, 
the corresponding stationary AF-algebras will be called   {\it companion AF-algebras}.

Let $N\ge 1$ be a natural number and  consider a (finite index) subgroup 
of the modular group given by the formula: 
\begin{equation}
\Gamma_1(N)=\left\{\left(\matrix{a & b\cr c & d}\right)\in SL_2({\Bbb Z}) ~|~
a\equiv d\equiv 1 ~\hbox{{\bf mod}} ~N,  ~c\equiv 0 ~\hbox{{\bf mod}}  ~N\right\}.
\end{equation}
Let ${\Bbb H}=\{z=x+iy\in {\Bbb C} ~|~ y>0\}$ be the upper half-plane  and 
let $\Gamma_1(N)$  act on ${\Bbb H}$  by the linear fractional
transformations;  consider an orbifold  ${\Bbb H}/\Gamma_1(N)$.
To compactify the orbifold 
at the cusps, one adds a boundary to ${\Bbb H}$,  so that 
${\Bbb H}^*={\Bbb H}\cup {\Bbb Q}\cup\{\infty\}$ and the compact Riemann surface 
$X_1(N)={\Bbb H}^*/\Gamma_1(N)$ is called a {\it modular curve}.   
The meromorphic functions $\phi(z)$ on ${\Bbb H}$ that
vanish at the cusps and such that
\begin{equation}
\phi\left({az+b\over cz+d}\right)= (cz+d)^2 \phi(z),\qquad
\forall \left(\matrix{a & b\cr c & d}\right)\in\Gamma_0(N), 
\end{equation}
are  called  {\it cusp forms} of weight two;  the (complex linear) space of such forms
will be denoted by $S_2(\Gamma_1(N))$.  The formula $\phi(z)\mapsto \omega=\phi(z)dz$ 
defines an isomorphism  $S_2(\Gamma_1(N))\cong \Omega_{hol}(X_1(N))$, where 
$\Omega_{hol}(X_1(N))$ is the space of all holomorphic differentials
on the Riemann surface $X_1(N)$.  Note that 
\linebreak
$\dim_{\Bbb C}(S_2(\Gamma_1(N))=\dim_{\Bbb C}(\Omega_{hol}(X_1(N))=g$,
where $g=g(N)$ is the genus of the surface $X_1(N)$. 
A {\it Hecke operator},  $T_n$, acts on $S_2(\Gamma_1(N))$ by the formula
$T_n \phi=\sum_{m\in {\Bbb Z}}\gamma(m)q^m$, where
$\gamma(m)= \sum_{a|\hbox{{\bf GCD}}(m,n)}a c_{mn/a^2}$ and 
$\phi(z)=\sum_{m\in {\Bbb Z}}c(m)q^m$ is the Fourier
series of the cusp form $\phi$ at $q=e^{2\pi iz}$.  Further,  $T_n$ is a
self-adjoint linear operator on the vector space $S_2(\Gamma_1(N))$
endowed with the Petersson inner product;  the algebra
${\Bbb T}_N :={\Bbb Z}[T_1,T_2,\dots]$ is a commutative algebra.
Any cusp form $\phi\in S_2(\Gamma_1(N))$ that is an eigenvector
for one (and hence all) of $T_n$, is referred  to
as a {\it Hecke eigenform}.
The Fourier  coefficients $c(m)$ of $\phi$ are algebraic integers,  and we 
denote by ${\cal K}_{\phi}={\Bbb Q}(c(m))$ an extension of the field ${\Bbb Q}$ 
by the Fourier coefficients of $\phi$.  Then ${\cal K}_{\phi}$
is a real algebraic number field of degree $1\le \deg~({\cal K}_{\phi} | {\Bbb Q})\le g$,
where $g$ is the genus of the surface $X_1(N)$, see e.g.  [Diamond \& Shurman 2005]   \cite{DS},   Proposition 6.6.4. 
Any embedding $\sigma: {\cal K}_{\phi}\to {\Bbb C}$ conjugates $\phi$ by acting
on its coefficients;  we write the corresponding Hecke  eigenform  
$\phi^{\sigma}(z):=\sum_{m\in {\Bbb Z}}\sigma(c(m))q^m$ and call $\phi^{\sigma}$
a {\it conjugate} of the Hecke eigenform $\phi$.

Let $\omega=\phi(z)dz\in\Omega_{hol}(X)$ be a holomorphic differential on a Riemann surface $X$. 
We shall  denote by $\Re~(\omega)$ a  closed  form on $X$  (the real part of $\omega$)
and  consider its periods $\lambda_i=\int_{\gamma_i}\Re~(\omega)$
against  a basis $\gamma_i$  in  the (relative) homology group 
$H_1(X,  Z(\Re~(\omega)); ~{\Bbb Z})$,  where $Z(\Re~(\omega))$ is the set of zeros of the form 
$\Re~(\omega)$.   
Assume $\lambda_i> 0$ and consider the vector $\theta=(\theta_1,\dots,\theta_{n-1})$
with $\theta_i=\lambda_{i+1} / \lambda_1$. The {\it Jacobi-Perron continued fraction} of
$\theta$ is given by the formula:
\begin{equation}
\left(\matrix{1\cr \theta}\right)=
\lim_{i\to\infty} \left(\matrix{0 & 1\cr I & b_1}\right)\dots
\left(\matrix{0 & 1\cr I & b_i}\right)
\left(\matrix{0\cr {\Bbb I}}\right)=
\lim_{i\to\infty} B_i\left(\matrix{0\cr {\Bbb I}}\right),
\end{equation}
where $b_i=(b^{(i)}_1,\dots, b^{(i)}_{n-1})^T$ is a vector of  non-negative integers,  
$I$ is the unit matrix and ${\Bbb I}=(0,\dots, 0, 1)^T$,   see e.g.  [Bernstein 1971]  \cite{BE}.
By ${\Bbb  A}_{\phi}$ we shall understand  the  AF-algebra
given the   Bratteli diagram  with  partial multiplicity matrices $B_i$. 
If $\phi(z)\in S_2(\Gamma_1(N))$ is a Hecke eigenform,  then the 
corresponding AF-algebra   ${\Bbb  A}_{\phi}$ is  {\it stationary}  with the partial multiplicity 
matrices $B_i=Const=B$;    moreover,  each   conjugate  eigenform  $\phi^{\sigma}$
defines  a  {\it companion} AF-algebra ${\Bbb A}_{\phi^{\sigma}}$.
 It is known that   $K_0^+({\Bbb A}_{\phi})\cong {\Bbb Z}+{\Bbb Z}\theta_1+\dots+{\Bbb Z}\theta_{n-1}
 \subset {\cal K}_{\phi}$,  where  ${\cal K}_{\phi}$  is an algebraic number field generated by the Fourier 
 coefficients of $\phi$,  see   \cite{Nik1}.

\section{Proof of theorem \ref{thm1}}
\begin{dfn}\label{dfn1}
Let $A_{\phi}\subset Jac~(X_1(fD))$ be an abelian variety associated to the
Hecke eigenform $\phi(z)\in S_2(\Gamma_1(fD))$,  see e.g. [Diamond \& Shurman 2005]  \cite{DS},  
Definition 6.6.3.   By $X_1^0(fD)$ we shall understand the Riemann surface of genus $g$,
such that 
\begin{equation}
Jac~(X_1^0(fD))\cong A_{\phi}.    
\end{equation}
By $\phi^0(z)dz\in\Omega_{hol}(X_1^0(fD))$ we denote the
image of the Hecke eigenform $\phi(z)dz\in\Omega_{hol}(X_1(fD))$
under the holomorphic map $X_1(fD)\to X_1^0(fD)$. 
\end{dfn}
\begin{rmk}
\textnormal{
The surface $X_1^0(fD)$ is correctly defined. Indeed, 
 since the abelian variety $A_{\phi}$ is the product of $g$ copies of an 
 elliptic curve with the complex multiplication, there exists a holomorphic map from 
 $A_{\phi}$ to the elliptic curve. For a Riemann
 surface $X$ of genus $g$ covering  the elliptic curve ${\cal E}_{CM}$ by a holomorphic
  map (such a surface and a map always exist),   one gets a period map $X\to A_{\phi}$
 by closing the arrows of a commutative diagram  $A_{\phi}\rightarrow {\cal E}_{CM}
 \leftarrow X$. It is easy to see,  that the Jacobian of $X$ coincides
 with  $A_{\phi}$ and we set  $X_1^0(fD) := X$.
}
\end{rmk}
\begin{lem}\label{lm1}
$g(X_1^0(fD))=\deg~({\cal K}^{ab}(k)~|~k)$.
\end{lem}
{\it Proof.}
By definition,  abelian variety $A_{\phi}$ is  the quotient of ${\Bbb C}^n$ by a lattice
of periods of the Hecke eigenform $\phi(z)\in S_2(\Gamma_1(fD))$  and all its
conjugates $\phi^{\sigma}(z)$ on the Riemann surface $X_1(fD)$.  These periods
are complex algebraic numbers generating the Hilbert class field ${\cal K}^{ab}(k)$
over imaginary quadratic field $k={\Bbb Q}(\sqrt{-D})$ modulo conductor $f$, see   
 [Hecke 1928]  \cite{Hec2},   [Shimura 1971] \cite{Shi2} and [Shimura 1972]  \cite{Shi1},   
 Section 8.  The number of linearly independent periods is equal to the total number
 of the conjugate eigenforms  $\phi^{\sigma}(z)$,  i.e. $|\sigma|=n=\dim_{\Bbb C} (A_{\phi})$.    
Since real dimension $\dim_{\Bbb R} (A_{\phi})=2n$,   we  conclude that 
$\deg~({\cal K}^{ab}(k)|{\Bbb Q})=2n$ and, therefore,   $\deg ~({\cal K}^{ab}(k)|k)=n$.
But $\dim_{\Bbb C} (A_{\phi})=g(X_1^0(fD))$ and one gets 
$g(X_1^0(fD))=\deg ~({\cal K}^{ab}(k)|k)$.  Lemma \ref{lm1} follows.
$\square$

\begin{cor}\label{cor1}
$g(X_1^0(fD))=|Cl~(R_f)|$.
\end{cor}
{\it Proof.}
Because  ${\cal K}^{ab}(k)$ is the Hilbert class field over  $k$ modulo conductor $f$,  
we must have 
\begin{equation}\label{eq12}
Gal~({\cal K}^{ab}(k)|k)\cong  Cl~(R_f),
\end{equation}
  where  $Gal$ is the Galois group of the extension ${\cal K}^{ab}(k)|k$  and  $Cl~(R_f)$
is the class group of ring  $R_f$,   see e.g.   [Silverman 1994]  \cite{S},  p.112.  
But $|Gal~({\cal K}^{ab}(k)|k)|=\deg ~({\cal K}^{ab}(k) | k)$ and  by lemma \ref{lm1} we have 
$\deg~({\cal K}^{ab}(k) | k)=g(X_1^0(fD))$.  
In view of this and isomorphism (\ref{eq12}),  one gets  $|Cl~(R_f)|=|Gal~({\cal K}^{ab}|k)|=g(X_1^0(fD))$.
Corollary  \ref{cor1} follows. 
$\square$

\begin{lem}\label{lm2}
$g(X_1^0(fD))=\deg~({\cal K}_{\phi} ~|~{\Bbb Q})$.
\end{lem}
{\it Proof.}
It is known that $\dim_{\Bbb C} (A_{\phi})=\deg~({\cal K}_{\phi} ~|~{\Bbb Q})$,
see e.g.  [Diamond \& Shurman 2005]  \cite{DS},   Proposition 6.6.4.  
But abelian variety $A_{\phi}\cong Jac~(X_1^0(fD))$ and, therefore,
$\dim_{\Bbb C} (A_{\phi})=\dim_{\Bbb C} (Jac~(X_1^0(fD)))=g(X_1^0(fD))$,
hence the lemma. 
$\square$

\begin{cor}\label{cor2}
$\deg~({\cal K}_{\phi} ~|~{\Bbb Q})=|Cl~({\goth R}_{\goth f})|$.
\end{cor}
{\it Proof.}
From lemma \ref{lm2} and corollary \ref{cor1} one gets $\deg~({\cal K}_{\phi} |{\Bbb Q})=|Cl~(R_f)|$.
In view of this and equality (\ref{eq3}),  one gets the conclusion of corollary \ref{cor2}.  
$\square$

\begin{lem}\label{lm3}
{\bf (Basic lemma)}
$Gal~({\cal K}_{\phi}~|~{\Bbb Q})\cong Cl~({\goth R}_{\goth f})$. 
\end{lem}
{\it Proof.}
Let us outline the proof.  In view of lemma \ref{lm2} and corollaries \ref{cor1}-\ref{cor2}, 
we  denote by $h$ the single integer $g(X_1^0(fD))=|Cl~(R_f)|=|Cl~({\goth R}_{\goth f})|=
\deg~({\cal K}_{\phi}|{\Bbb Q})$.   
Since $\deg~({\cal K}_{\phi}|{\Bbb Q})=h$,  there exist $\{\phi_1,\dots,\phi_h\}$ conjugate
Hecke eigenforms $\phi_i(z)\in S_2(\Gamma_1(fD))$,   see e.g.  
[Diamond \& Shurman 2005]  \cite{DS}, Theorem 6.5.4;    thus  one gets $h$
holomorphic forms  $\{\phi_1^0,\dots,\phi_h^0\}$ on the Riemann surface $X_1^0(fD)$. 
Let    $\{{\Bbb A}_{\phi_1^0},\dots, {\Bbb A}_{\phi_h^0}\}$
be  the corresponding stationary AF-algebras;   the ${\Bbb A}_{\phi_i^0}$ are {\it companion} 
AF-algebras,  see Section 1.2.  
Recall that the characteristic polynomial for the  partial multiplicity matrices $B_{\phi_i^0}$
of  companion AF-algebras ${\Bbb A}_{\phi_i^0}$  is the same;  it is a minimal polynomial
of degree $h$ and  let $\{\lambda_1,\dots,\lambda_h\}$ be the roots of such a polynomial,
compare with  [Effros 1981]  \cite{E},  Corollary 6.3.
Since $\det~(B_{\phi_i^0})=1$,  the numbers $\lambda_i$ are algebraic units of the field ${\cal K}_{\phi}$.   
Moreover,  $\lambda_i$ are algebraically conjugate and can be taken for generators
of the extension ${\cal K}_{\phi}|{\Bbb Q}$;   since $\deg~({\cal K}_{\phi}|{\Bbb Q})=h=|Cl~({\goth R}_{\goth f})|$
there exists a natural action of group $Cl~({\goth R}_{\goth f})$ on these generators.   The action extends  to
automorphisms  of the entire field    ${\cal K}_{\phi}$ preserving ${\Bbb Q}$;  thus one gets
the Galois group of   extension ${\cal K}_{\phi}|{\Bbb Q}$ and  an isomorphism  
$Gal~({\cal K}_{\phi}|{\Bbb Q})\cong Cl~({\goth R}_{\goth f})$.  Let us pass to a step-by-step argument.

\bigskip
(i)    Let   $h:=g(X_1^0(fD))=|Cl~(R_f)|=|Cl~({\goth R}_{\goth f})|$
 and let $\phi(z)\in S_2(\Gamma_1(fD))$ be the 
Hecke eigenform.  It is known that there exists  $\{\phi_1,\dots,\phi_h\}$ 
conjugate Hecke eigenforms,  so that $\phi(z)$ is one of them,  see  
[Diamond \& Shurman 2005]  \cite{DS},   Theorem 6.5.4.
Let $\{\phi_1^0,\dots,\phi_h^0\}$  be the corresponding forms on the Riemann surface $X_1^0(fD)$.  
\begin{rmk}
\textnormal{
The forms  $\{\phi_1^0,\dots,\phi_h^0\}$ can be taken for a basis in the space $\Omega_{hol}(X_1^0(fD))$;
we leave it to the reader to verify,  that abelian variety $A_{\phi}$ is isomorphic to the quotient of ${\Bbb C}^h$ by
the lattice  of  periods of  holomorphic differentials $\phi_i^0(z) dz$ on $X_1^0(fD)$. 
}
\end{rmk}

\bigskip
(ii)  Let ${\Bbb A}_{\phi_i^0}$ be the AF-algebra corresponding to holomorphic differential $\phi_i^0(z)dz$
on $X_1^0(fD)$,  see Section 2.2;
the set  $\{{\Bbb A}_{\phi_1^0},\dots, {\Bbb A}_{\phi_h^0}\}$   consists of the companion AF-algebras.
It is known  that each  ${\Bbb A}_{\phi_i^0}$ is a stationary AF-algebra,  i.e. its  partial multiplicity matrix is a constant;
we shall denote such a matrix  by $B_{\phi_i^0}$.

\bigskip
(iii) By definition,  the  matrices $B_{\phi_i^0}$ of  companion $AF$-algebras  ${\Bbb A}_{\phi_i}$
have the same characteristic polynomial $p(x)\in {\Bbb Z}[x]$;   the matrices $B_{\phi_i^0}$ itself are not pairwise
similar and,  therefore,   the AF-algebras   ${\Bbb A}_{\phi_i^0}$ are not pairwise isomorphic.  The total
number $h$ of such matrices is equal to the class number of the endomorphism ring of 
pseudo-lattice $K_0^+({\Bbb A}_{\phi_i^0})\cong {\Bbb Z}+{\Bbb Z}\theta_1^i+\dots+{\Bbb Z}\theta_{h-1}^i
\subset {\cal K}_{\phi}$,   see  [Effros 1981]  \cite{E},  Corollary 6.3.   
\begin{rmk}\label{rmk4}
\textnormal{
Notice that there are $\{X_1,\dots, X_h\}$ pairwise non-isomorphic Riemann surfaces
$X:=X_1^0(fD)$ endowed with a holomorphic map $X_i\to {\cal E}_i$,
where $\{{\cal E}_1,\dots, {\cal E}_h\}$ are pairwise non-isomorphic elliptic curves 
${\cal E}_{CM}^{(-D,f)}$ corresponding to elements of the group $Cl~(R_f)$.
Thus the companion AF-algebras $\{{\Bbb A}_{\phi_1^0},\dots, {\Bbb A}_{\phi_h^0}\}$ 
can be viewed as  coordinate rings of  $\{X_1,\dots, X_h\}$;
the latter means that   ${\Bbb A}_{\phi_i^0}$ discern non-isomorphic Riemann  surfaces and 
$K_0^+({\Bbb A}_{\phi_i^0})\cong {\Bbb Z}+{\Bbb Z}\theta_1^i+\dots+{\Bbb Z}\theta_{h-1}^i$
represents the moduli space of $X_1^0(fD)$.   
}
\end{rmk}

\bigskip
(iv)  The polynomial $p(x)$ is minimal and splits in the totally real  field ${\cal K}_{\phi}$.  
Indeed,  matrices $B_{\phi_i^0}$ generate the Hecke algebra ${\Bbb T}_N$
on $S_2(\Gamma_1(N))$;  thus each $B_{\phi_i^0}$ is self-adjoint and,  therefore,
all eigenvalues are real of multiplicity one;   since $B_{\phi_i^0}$ is integer, all roots
of characteristic polynomial $p(x)$ of    $B_{\phi_i^0}$ belong to the field ${\cal K}_{\phi}$.

\bigskip
(v)  Let  $p(x)=(x-\lambda_1)\dots (x-\lambda_h)$.    It is easy to see  that  $\lambda_i$ are algebraic units
of the field ${\cal K}_{\phi}$ because $\det~(B_{\phi_i^0})=1$;   note  that  numbers 
 $\{\lambda_1,\dots,\lambda_h\}$ are  algebraically conjugate.   
 Since $\deg~({\cal K}_{\phi}|{\Bbb Q})=h$,   the numbers $\lambda_i$  can be taken for
generators of the field ${\cal K}_{\phi}$,   i.e.    ${\cal K}_{\phi}={\Bbb Q}(\lambda_1,\dots,\lambda_h)$.

\bigskip
(vi) Finally,  let us  establish  an explicit formula for the isomorphism 
\begin{equation}\label{eq17}
Cl~({\goth R}_{\goth f})\to Gal~({\cal K}_{\phi}|{\Bbb Q})
\end{equation}
Since $Gal~({\cal K}_{\phi}|{\Bbb Q})$  is an automorphism group of the
field ${\cal K}_{\phi}$ preserving ${\Bbb Q}$,   it will suffice  to define  the action $\ast$  of an 
element $a\in Cl~({\goth R}_{\goth f})$  on  the generators $\lambda_i$ of ${\cal K}_{\phi}$.   
Let $\{a_1,\dots,a_h\}$ be the set of all elements of the group $Cl~({\goth R}_{\goth f})$.
For an element $a\in Cl~({\goth R}_{\goth f})$ define an index function $\alpha$ 
by the formula $a_ia=a_{\alpha(i)}$.   Then the action $\ast$  of an 
element $a\in Cl~({\goth R}_{\goth f})$ on the generators $\lambda_i$ of the field
${\cal K}_{\phi}$ is given by the formula:
\begin{equation}\label{eq18}
a\ast \lambda_i:=  \lambda_{a(i)},  \qquad\forall a\in Cl~({\goth R}_{\goth f}).
\end{equation}
It is easy  to verify that formula (\ref{eq18})  gives an isomorphism
$Cl~({\goth R}_{\goth f})\to Gal~({\cal K}_{\phi}|{\Bbb Q})$,  which is  independ of the choice of
$\{a_i\}$ and $\{\lambda_i\}$.  This argument completes  the proof of lemma \ref{lm3}.
$\square$

\begin{rmk}\label{rmk5}
\textnormal{
The class field theory says that ${\goth f}=f^m$, i.e. 
the extensions of fields  $k$ and ${\goth k}$ must ramify over the same 
set of prime ideals.  Indeed, consider the commutative diagram below, 
\begin{figure}[here]
\begin{picture}(300,110)(-100,-5)
\put(30,70){\vector(0,-1){35}}
\put(155,70){\vector(0,-1){35}}
\put(52,23){\vector(1,0){73}}
\put(52,83){\vector(1,0){73}}
\put(28,20){$I_{\goth f}$}
\put(135,20){$Gal~({\cal K}_{\phi}|{\Bbb Q})$}
\put(28,80){$I_f$}
\put(135,80){$Gal~({\cal K}^{ab}(k)|{\Bbb Q})$}
\put(68, 90){Artin}
\put(50, 70){homomorphism}
\put(68, 30){Artin}
\put(50, 10){homomorphism}
\end{picture}
\end{figure}
where $I_f$ and $I_{\goth f}$ are groups of all ideals of $k$ and
${\goth k}$,  which are relatively prime to the principal ideals $(f)$ and $({\goth f})$,
respectively. Since $Gal~({\cal K}^{ab}(k)|{\Bbb Q})\cong Gal~({\cal K}_{\phi}|{\Bbb Q})$, 
one gets an isomorphism $I_f\cong I_{\goth f}$,  i.e. ${\goth f}=f^m$ for some
positive integer $m$.   
}
\end{rmk}
\begin{cor}\label{cr3}
The Hilbert class field of real quadratic field ${\goth k}={\Bbb Q}(\sqrt{D})$
modulo conductor ${\goth f}\ge 1$ is isomorphic to the field ${\goth k}({\cal K}_{\phi})$ 
generated by the Fourier coefficients of the Hecke eigenform $\phi(z)\in S_2(\Gamma_1(fD))$.
\end{cor}
{\it Proof.}
As in the classical case of imaginary quadratic fields,   
notice that    $\deg~({\cal K}_{\phi}|{\Bbb Q})=\deg~({\goth k}({\cal K}_{\phi})|{\goth k})=Cl~({\goth R}_{\goth f})$;
therefore  corollary \ref{cr3} is an implication of lemma \ref{lm3} and isomorphism  
$Gal~({\cal K}_{\phi}|{\Bbb Q})\cong  Gal~({\goth k}({\cal K}_{\phi})|{\goth k})\cong Cl~({\goth R}_{\goth f})$. 
$\square$

\bigskip
Theorem \ref{thm1} follows from corollary \ref{cr3}.
$\square$

\section{Examples}
Along with the method of Stark's units [Cohen \& Roblot 2000] \cite{CoRo1},  theorem \ref{thm1} can be used in the  computational number theory.  For the sake of clarity,  
we shall  consider  the simplest   examples;  the rest can be found in Figure 1.  
\begin{exm}
\textnormal{
Let $D=15$.    The class number of quadratic field $k={\Bbb Q}(\sqrt{-15})$ 
is known to be  $2$;   such a  number for quadratic field ${\goth k}={\Bbb Q}(\sqrt{15})$
is also equal to $2$. Thus,
\begin{equation}
Cl~({\goth R}_{{\goth f}=1})\cong Cl~(R_{f=1})\cong {\Bbb Z}/2 {\Bbb Z},  
\end{equation}
and isomorphism (\ref{eq3}) is trivially satisfied for each power $m$,
i.e. one obtains an unramified extension. 
By theorem \ref{thm1},   the Hilbert class field of ${\goth k}$
 is generated by the Fourier coefficients of the Hecke eigenform
 $\phi(z)\in S_2(\Gamma_1(15))$. 
  Using the computer program {\it SAGE}  created by  William ~A. ~Stein,  one finds an
  irreducible factor $p(x)=x^2-4x+5$ of the  characteristic polynomial of the Hecke operator $T_{p=2}$ 
acting on the space    $S_2(\Gamma_1(15))$.  Therefore, the Fourier coefficient $c(2)$   coincides  with a root 
 of equation  $p(x)=0$;   in other words,   we arrive at an extension of ${\goth k}$ by the polynomial  $p(x)$.
 The  generator $x$ of the field ${\cal K}_{\phi}={\Bbb Q}(c(2))$ is a root of the  
 bi-quadratic equation $[(x-2)^2+1]^2-15=0$;   it is easy to see that  
 $x=2+\sqrt{-1+\sqrt{15}}$.  
 One  concludes,   that the field 
${\cal K}_{\phi}\cong {\Bbb Q}\left(\sqrt{-1+\sqrt{15}}\right)$
is the Hilbert class field of quadratic field   ${\goth k}={\Bbb Q}(\sqrt{15})$.    
 }
\end{exm}
\begin{exm}
\textnormal{
Let $D=14$.   It is known,  that for the quadratic field $k={\Bbb Q}(\sqrt{-14})$ we have $|Cl~(R_{f=1})|=4$,
while for the  quadratic field ${\goth k}={\Bbb Q}(\sqrt{14})$ it holds $|Cl~({\goth R}_{{\goth f}=1})|=1$.
However, for the ramified extensions one obtains
the following isomorphism:  
\begin{equation}
Cl~({\goth R}_{{\goth f}=2^3})\cong Cl~(R_{f=2})\cong  {\Bbb Z}/4 {\Bbb Z},  
\end{equation}
where $m=3$ is the smallest integer satisfying formula (\ref{eq3}).  
 By theorem \ref{thm1},  the Hilbert class field of ${\goth k}$  modulo ${\goth f}=8$ is generated by the 
 Fourier coefficients of the Hecke eigenform  $\phi(z)\in S_2(\Gamma_1(2\times 14))$.  
   Using the  {\it SAGE},  one finds that  the  characteristic polynomial of the Hecke operator $T_{p=3}$ on   
 $S_2(\Gamma_1(2\times 14))$  has an irreducible factor $p(x)=x^4+3x^2+9$.
 Thus the Fourier coefficient $c(3)$  is a root of the polynomial $p(x)$ and one gets 
 an extension of ${\goth k}$ by the polynomial  $p(x)$.
 In other words,  generator $x$ of the field ${\cal K}_{\phi}={\Bbb Q}(c(3))$ 
 is a root of the polynomial equation $(x^4+3x^2+9)^2-4\times 14=0$.
 The bi-quadratic equation $x^4+3x^2+9-2\sqrt{14}=0$ has discriminant $-27+8\sqrt{14}$
 and one finds a generator of ${\cal K}_{\phi}$ to be  $^4\sqrt{-27+8\sqrt{14}}$. 
 Thus  the field  ${\Bbb Q}\left(^4\sqrt{-27+8\sqrt{14}}\right)$
is the Hilbert class field over  ${\Bbb Q}(\sqrt{14})$    modulo  conductor ${\goth f}=8$.
Clearly, the extension is ramified over the prime ideal ${\goth p}=(2)$.      
 }
\end{exm}
\begin{rmk}
\textnormal{
Table 1 below lists   quadratic fields for some  square-free discriminants $2\le D\le 101$. 
 The conductors $f$ and ${\goth f}$ satisfying equation (\ref{eq3}) were calculated using 
 tables for  the  class number  of   non-maximal orders in  quadratic fields  
 posted  at {\sf www.numbertheory.org};  the  site is   maintained by Keith ~Matthews.
 We focused on  small conductors;   the interested reader 
can  compute  the higher conductors  using  a pocket  calculator.    
In contrast,    computation  of generator $x$   of the Hilbert   class field
 require  the online program {\it SAGE} created    
by William ~A. ~Stein.   We write  an explicit
formula for  $x$ or its minimal polynomial $p(x)$ over ${\goth k}$.  
}
\end{rmk}

\begin{figure}[top]
\begin{tabular}{c|c|c|c|c}
\hline
&&&&\\
$D$ & $f$ & $Cl~(R_f)$ & ${\goth f}$ &  Hilbert class field of ${\Bbb Q}(\sqrt{D})$\\
&&&& modulo conductor ${\goth f}$\\
&&&&\\
\hline
$2$ & $1$ & trivial & $1$ & ${\Bbb Q}(\sqrt{2})$ \\
\hline
$3$ & $1$ & trivial & $1$ &  ${\Bbb Q}(\sqrt{3})$ \\
\hline
$7$ & $1$ & trivial  & $1$ & ${\Bbb Q}(\sqrt{7})$ \\
\hline
$11$ & $1$ & trivial & $1$ & ${\Bbb Q}(\sqrt{11})$ \\
\hline
&&&&\\
$14$ & $2$ & ${\Bbb Z}/4 {\Bbb Z}$ & $8$ &  ${\Bbb Q}\left(^4\sqrt{-27+8\sqrt{14}}\right)$ \\
&&&&\\
\hline
&&&&\\
$15$ & $1$ & ${\Bbb Z}/2 {\Bbb Z}$ & $1$ &   ${\Bbb Q}\left(\sqrt{-1 +\sqrt{15}}\right)$ \\
&&&&\\
\hline
$19$ & $1$ & trivial & $1$ & ${\Bbb Q}(\sqrt{19})$\\
\hline
$21$ & $2$ & ${\Bbb Z}/4 {\Bbb Z}$ & $8$ & ${\Bbb Q}\left(^4\sqrt{-3 +2\sqrt{21}}\right)$ \\
\hline
&&&&\\
$35$ & $1$ & ${\Bbb Z}/2 {\Bbb Z}$ & $1$ & ${\Bbb Q}\left(\sqrt{17+\sqrt{35}}\right)$ \\
&&&&\\
\hline
$43$ & $1$ & trivial & $1$ & ${\Bbb Q}(\sqrt{43})$\\
\hline
&&&&\\
$51$ & $1$ & ${\Bbb Z}/2 {\Bbb Z}$ & $1$ & ${\Bbb Q}\left(\sqrt{17+\sqrt{51}}\right)$ \\
&&&&\\
\hline
$58$ & $1$ & ${\Bbb Z}/2 {\Bbb Z}$ & $1$ & ${\Bbb Q}\left(\sqrt{-1+\sqrt{58}}\right)$ \\
\hline
$67$  & $1$ & trivial & $1$ & ${\Bbb Q}(\sqrt{67})$ \\
\hline
$82$ & $1$ & ${\Bbb Z}/4 {\Bbb Z}$ & $1$ & $x^4-2x^3+4x^2-8x+16$ \\
\hline
&&&&\\
$91$ & $1$ & ${\Bbb Z}/2 {\Bbb Z}$ & $1$ & ${\Bbb Q}\left(\sqrt{-3 +\sqrt{91}}\right)$ \\
&&&&\\
\hline
\end{tabular}
\caption{Square-free discriminants $2\le D\le 101$.}
\end{figure}

\bigskip\noindent
{\sf Acknowledgment.}  I   thank  Yu.~I.~Manin   for helpful correspondence.




\vskip1cm

\textsc{The Fields Institute for Research in Mathematical Sciences, Toronto, ON, Canada,  
E-mail:} {\sf igor.v.nikolaev@gmail.com}


\end{document}